# A NEW NUMERICAL TECHNIQUE FOR SOLVING FRACTIONAL PARTIAL DIFFERENTIAL EQUATIONS

**Omer Acan[1*], Dumitru Baleanu[2,3]**

[1] *Department of Mathematics, Faculty of Art and Science, Siirt University, 56100, Siirt, Turkey*
[2] *Department of Mathematics and Computer Sciences, Faculty of Art and Science, Balgat 06530, Ankara, Turkey*
[3] *Institute of Space Science, Magurele-Bucharest, Romania*

**Abstract**

we propose conformable Adomian decomposition method (CADM) for fractional partial differential equations (FPDEs) at the first time. This method is a new Adomian decomposition method based on conformable derivative to solve FPDEs. At the same time, conformable reduced differential transform method (CRDTM) for FPDEs is briefly given and a numerical comparison is made between this method and the newly introduced CADM. In applied science, CADM can be used as an alternative method to obtain approximate and analytical solutions for FPDEs as CRDTM. In this study, linear and non-linear three problems are solved by these two methods. In these schemes, the obtained solutions take the form of a convergent series with easily computable algorithms. For the applications, the obtained results by these methods are compared to each other and with the exact solutions. When applied to FPDEs, it is seem that CADM approach produces easy, fast and reliable solutions as CFRDTM.

## 1. Introduction

Linear and non-linear both fractional and non-fractional problems of differential equations play a major role in various fields such as biology, physics, chemistry, mathematics, astronomy, fluids mechanics, optics, applied mathematics, and engineering. It is not always possible to find analytical solutions to these problems [1–10]. Therefore, it is very important to handle these problems appropriately and solve them or develop solutions. Therefore, it is very important to handle these problems appropriately and solve them or develop solutions. Adomian decomposition method (ADM), which is introduced [11–13] in the 1980's, is one of the important mathematical methods used to solve many problems in real world. Since then, Since then, a number of studies have been conducted on ADM such as linear and nonlinear, homogeneous and nonhomogeneous operator equations, including fractional or nonfractional ordinary differential equations (ODEs), partial differential equations (PDEs), integral equations, integro-differential equations, etc. (see [14–24] and references therein). Recently, a new derivative called CDO was introduced [25–27] In science, by the help of this new derivative, the behaviors of many problems has been studied and some solutions techniques have been developed [27–37]. This new subject gives academicians an opportunity to study further in many engineering, physical and applied mathematics problems.

The aim of this study is to introduce CADM by using CDO and ADM for the first time in the literature. This method can be used to solve many linear and non-linear FPDEs. We will briefly mentioned CRDTM to compare our CADM with it. The problems will be solved both by the CRDTM and the first proposed CADM. The obtained solutions by these methods will be compared. For this, in section 2, we give some basic definitions and important properties of CDO. In section 3, we propose CADM. In sections 4, we introduce CRDTM to compare with our method. In section 5, we give applications of CADM and CRDTM. We give the conclusion in the final section.

## 2. On The Conformable Derivative

**Definition 2.1.**
Given a functio $f:[0,\infty) \to \mathbb{R}$. Then the conformable derivative of $f$ order $\alpha$ is defined by [25,26]

$$(T_\alpha f)(x) = \lim_{\varepsilon \to 0} \frac{f(x + \varepsilon x^{1-\alpha}) - f(x)}{\varepsilon}$$

for all $x > 0$, $\alpha \in (0,1]$.

---

[*] Corresponding author.
E-mail addresses: omeracan@yahoo.com (O. Acan), dumitru@cankaya.edu.tr (D. Baleanu)

**Lemma 2.1.** [25,26]

Let $\alpha \in (0,1]$ and $f, g$ be $\alpha$-differentiable at a point $x > 0$. Then

(i) $T_\alpha(af + bg) = a(T_\alpha f) + b(T_\alpha g)$ for $a, b \in \mathbb{R}$,

(ii) $T_\alpha(x^p) = px^{p-\alpha}$, for all $p \in \mathbb{R}$,

(iii) $T_\alpha(f(x)) = 0$, for all constant functions $f(x) = \lambda$,

(iv) $T_\alpha(fg) = f(T_\alpha g) + g(T_\alpha f)$,

(v) $T_\alpha(f/g) = \dfrac{g(T_\alpha f) - f(T_\alpha g)}{g^2}$,

(vi) If $f(x)$ is differentiable, then $T_\alpha(f(x)) = x^{1-\alpha} \dfrac{d}{dx} f(x)$.

**Definition 2.2.**

Given a function $f : [a, \infty) \to \mathbb{R}$. Then the conformable derivative of $f$ order $\alpha$ is defined by [25,26]:

$$(T_\alpha^a f)(x) = \lim_{\varepsilon \to 0} \frac{f(x + \varepsilon(x-a)^{1-\alpha}) - f(x)}{\varepsilon}$$

for all $x > 0$, $\alpha \in (0,1]$.

**Definition 2.3.**

Let $f$ be an $n$ times differentiable at $x$. Then the conformable derivative of $f$ order $\alpha$ is defined by [25,26]:

$$(T_\alpha f)(x) = \lim_{\varepsilon \to 0} \frac{f^{(\lceil \alpha \rceil - 1)}(x + \varepsilon x^{(\lceil \alpha \rceil - \alpha)}) - f^{(\lceil \alpha \rceil - 1)}(x)}{\varepsilon}$$

for all $x > 0$, $\alpha \in (n, n+1]$, $\lceil \alpha \rceil$ is the smallest integer greater than or equal to $\alpha$.

**Lemma 2.2.** [25,26]

Let $f$ be an $n$ times differentiable at $x$. Then

$$T_\alpha(f(x)) = x^{\lceil \alpha \rceil - \alpha} f^{(\lceil \alpha \rceil)}(x)$$

for all $x > 0$, $\alpha \in (n, n+1]$.

### 3. Conformable Adomian Decomposition Method

We will briefly introduce *CADM* for non-linear *FPDEs* in this section. For the purpose of illustration of the methodology to the proposed methods, we write the non-linear *FPDEs* in the standard operator form

$$L_\alpha(u(x,t)) + R(u(x,t)) + N(u(x,t)) = g(x,t) \tag{3.1}$$

where $L_\alpha = {}_\alpha T$ is a linear operator with conformable derivative of order $\alpha$ ($n < \alpha \leq n+1$), $R$ is the other part of the linear operator, $N$ is a non-linear operator and $g(x,t)$ is a non-homogeneous term. If the linear operator in eq. (3.1) is applied to Lemma 2.2, the following equation is obtained:

$$t^{\lceil \alpha \rceil - \alpha} \frac{\partial^{\lceil \alpha \rceil}}{\partial t^{\lceil \alpha \rceil}} u(x,t) + R(u(x,t)) N(u(x,t)) = g(x,t) \tag{3.2}$$

Applying $L_\alpha^{-1} = \int_0^t \int_0^{\xi_1} \cdots \int_n^{\xi_{n-1}} \frac{1}{\xi_n^{\lceil \alpha \rceil - \alpha}} (.) d\xi_n d\xi_{n-1} \cdots d\xi_1$, ($n < \alpha \leq n+1$) the inverse of operator, to both sides of (3.2), it is obtained as

$$L_\alpha^{-1} L_\alpha(u(x,t)) = L_\alpha^{-1} g(x,t) - L_\alpha^{-1} R(u(x,t)) - L_\alpha^{-1} N(u(x,t)). \tag{3.3}$$

The general solution of the given equation is decomposed into the sum

$$u(x,t) = \sum_{n=0}^{\infty} u_n(x,t). \tag{3.4}$$

The nonlinear term $N(u)$ can be decomposed into the infinite series of polynomial obtained by

$$N(u) = \sum_{n=0}^{\infty} A_n, \ (u_0, u_1, \ldots, u_n), \tag{3.5}$$

where $A_n$ is the so-called Adomian polynomials. These polynomials can be calculated for all forms of nonlinearity according to specific algorithms constructed by Adomian [12,17]. $u$ and $N(u)$, respectively, is obtained as

$$u = \sum_{i=0}^{\infty} \lambda^i u_i, \quad N(u) = N\left(\sum_{i=0}^{\infty} \lambda^i u_i\right) = \sum_{i=0}^{\infty} \lambda^i A_i \tag{3.6}$$

where $\lambda$ is the convenience parameter. From (3.6), Adomian polynomials $A_n$ are obtained as

$$n! A_n = \frac{d^n}{d\lambda^n}\left[N\left(\sum_{n=0}^{\infty} \lambda^n u_n\right)\right]_{\lambda=0}.$$

These polynomials can be calculated easily with the following Maple code:

```
restart:
ord:=m:    #Order (Enter a number for m.)
NF:=N(u(x,t))  #Nonlinear Function
u[t]:=sum(u[b]*t^b,b=0..ord):
NF[t]:=subs(u(x,t)=u[t],NF):
s:=expand(NF[t],t):
dt:=unapply(s,t):
for i from 0 to ord do
A[i]:=((D@@i)(dt)(0)/i!):
print(AA[i],A[i]): #Adomian Polinoms
od:
```

Substituting (3.4) and (3.5) into (3.3), it is obtained

$$\sum_{n=0}^{\infty} u_n = \theta + L_\alpha^{-1} g - L_\alpha^{-1} R\left(\sum_{n=0}^{\infty} u_n\right) - L_\alpha^{-1}\left(\sum_{n=0}^{\infty} A_n\right). \tag{3.7}$$

where $\theta = u(x,0)$ is initial condition. From (3.7), the iterates are defined by the following recursive way:

$$\begin{aligned} u_0 &= \theta + L_\alpha^{-1} g, \\ u_1 &= -L_\alpha^{-1} R\, u_0 - L_\alpha^{-1} A_0, \\ &\vdots \\ u_{n+1} &= -L_\alpha^{-1} R\, u_n - L_\alpha^{-1} A_n, \quad n \geq 0. \end{aligned} \tag{3.8}$$

Therefore, the approximation solution of (3.1) is obtained by

$$\tilde{u}_m(x,t) = \sum_{n=0}^{m} u_n(x,t).$$

Hence the exact solution of (3.1) given as

$$u(x,t) = \lim_{m \to \infty} \tilde{u}_m(x,t).$$

## 4. Conformable Reduced Differential Transform Method

In this section, it is given basic definitions and properties of CFRDTM for fractional PDEs. The lowercase $u(x,t)$ represent the original function while the uppercase $U_k^\alpha(x)$ stand for the conformable fractional reduced differential transformed (CRDT) function.

**Definition 4.1.** [36]

Assume $u(x,t)$ is analytic and differentiated continuously with respect to time $t$ and space $x$ in the its domain. CFRDT of $u(x,t)$ is defined as

$$U_k^\alpha(x) = \frac{1}{\alpha^k k!}\left[\left(T_\alpha^{(k)} u\right)\right]_{t=t_0}$$

where some $0 < \alpha \leq 1$, $\alpha$ is a parameter describing the order of conformable fractional derivative, $T_\alpha^{(k)} u = \underbrace{(T_\alpha\, T_\alpha \cdots T_\alpha)}_{k \text{ times}} u(x,t)$ and the $t$ dimensional spectrum function $U_k^\alpha(x)$ is the CFRDT function.

**Definition 4.2.** [36]

Let $U_k^\alpha(x)$ be the CFRDT of $u(x,t)$. Inverse CFRDT of $U_k^\alpha(x)$ is defined as

$$u(x,t) = \sum_{k=0}^{\infty} U_k^\alpha(x)(t-t_0)^{\alpha k} = \sum_{k=0}^{\infty} \frac{1}{\alpha^k k!}\left[{}_tT_\alpha^{(k)}u\right]_{t=t_0}(t-t_0)^{\alpha k}.$$

CFRDT of initial conditions for integer order derivatives are defined as

$$U_k^\alpha(x) = \begin{cases} \dfrac{1}{(\alpha k)!}\left[\dfrac{\partial^{\alpha k}}{\partial t^{\alpha k}}u(x,t)\right]_{t=t_0} & if\ \alpha k \in \mathbb{Z}^+ \\ 0 & if\ \alpha k \notin \mathbb{Z}^+ \end{cases} \quad for\ k=0,1,2,\ldots,\left(\dfrac{n}{\alpha}-1\right)$$

where $n$ is the order of conformable fractional PDE.

By consideration of $U_0^\alpha(x) = f(x)$ as transformation of initial condition

$$u(x,0) = f(x).$$

A straightforward iterative calculations gives the $U_k^\alpha(x)$ values for $k=1,2,3,\ldots,n$. Then the inverse transformation of the $\{U_k^\alpha(x)\}_{k=0}^n$ gives the approximation solution as:

$$\tilde{u}_n(x,t) = \sum_{k=0}^{n} U_k^\alpha(x) t^{k\alpha},$$

where n is order of approximation solution. The exact solution is given by:

$$u(x,t) = \lim_{n\to\infty} \tilde{u}_n(x,t)$$

The fundamental operations of CFRDT that can be deduced from Definition 4.1 and Definition 4.2 are listed Table 4.1.

| Table 4.1. Basic operations of CFRDTM [36]. | |
|---|---|
| **Original function** | **Transformed function** |
| $u(x,t)$ | $U_k^\alpha(x) = \dfrac{1}{\alpha^k k!}\left[\left({}_tT_\alpha^{(k)}u\right)\right]_{t=t_0}$ |
| $u(x,t) = av(x,t) \pm bw(x,t)$ | $U_k^\alpha(x) = aV_k^\alpha(x) \pm bW_k^\alpha(x)$ |
| $u(x,t) = v(x,t)w(x,t)$ | $U_k^\alpha(x) = \sum_{s=0}^{k} V_s^\alpha(x) W_{k-s}^\alpha(x)$ |
| $u(x,t) = {}_tT_\alpha v(x,t)$ | $U_k^\alpha(x) = \alpha(k+1)V_{k+1}^\alpha(x)$ |
| $u(x,t) = x^m(t-t_0)^n$ | $U_k^\alpha(x) = x^m \delta\left(k - \dfrac{n}{\alpha}\right)$ |

## 5. Applications

To illustrate the effectiveness of the given CADM and CRDTM, three examples are considered in this section. All the results are calculated by software MAPLE.

**Example 4.1:**

Firstly, Consider the linear time and space fractional diffusion equation:

$$\frac{\partial^\alpha}{\partial t^\alpha}u(x,t) = \frac{\partial^{2\beta}}{\partial x^{2\beta}}u(x,t) \quad t>0, x \in R, 0<\alpha,\beta \leq 1 \tag{5.1}$$

subject to the initial condition

$$u(x,0) = \sin\left(\frac{x^\beta}{\beta}\right). \tag{5.2}$$

Exact solution of the problem) in conformable sense is

$$u(x,t) = \sin\left(\frac{x^\beta}{\beta}\right)e^{-\frac{t^\alpha}{\alpha}}.$$

***Solution by CADM:***

Solve this problem by using CADM. Let $L_\alpha = T_\alpha = \frac{\partial^\alpha}{\partial t^\alpha}$ be a linear operator, then the operator form of (5.1) is as follows

$$T_\alpha u(x,t) = \frac{\partial^{2\beta}u(x,t)}{\partial x^{2\beta}} \quad t>0, x \in R, 0<\alpha,\beta \leq 1 \tag{5.3}$$

By the help of Lemma 2.2, eq. (4.4) can be written as

$$t^{1-\beta}\frac{\partial u(x,t)}{\partial t} = \frac{\partial^{2\beta}u(x,t)}{\partial x^{2\beta}} \quad t>0, x \in R, 0<\alpha,\beta \leq 1. \tag{5.4}$$

If $L_\alpha^{-1} = \int_0^t \frac{1}{\xi^{1-\alpha}}(.)d\xi$, which is the inverse of $L_\alpha$, is applied to both sides of eq. (5.4), we get

$$u(x,t) = u(x,0) - L_\alpha^{-1}\left(\frac{\partial^{2\beta}}{\partial x^{2\beta}}u(x,t)\right).$$

According to (3.8) and the initial condition (5.2), we can write

$$\begin{aligned} u_0 &= \sin\left(\frac{x^\beta}{\beta}\right), \\ u_1 &= -L_\alpha^{-1}\left(\frac{\partial^2}{\partial x^2}u_0\right) \\ &\vdots \\ u_{n+1} &= -L_\alpha^{-1}\left(\frac{\partial^2}{\partial x^2}u_n\right), n \geq 0. \end{aligned} \tag{5.5}$$

From (5.5), we conclude the terms of decomposition series as:

$$\begin{aligned} u_0 &= \sin\left(\frac{x^\beta}{\beta}\right), \ u_1 = -\sin\left(\frac{x^\beta}{\beta}\right)\frac{t^\alpha}{\alpha}, \ u_2 = \sin\left(\frac{x^\beta}{\beta}\right)\frac{t^{2\alpha}}{2\alpha^2}, \\ u_3 &= -\sin\left(\frac{x^\beta}{\beta}\right)\frac{t^{3\alpha}}{3!\alpha^3}, \cdots, u_n = (-1)^n \sin\left(\frac{x^\beta}{\beta}\right)\frac{t^{n\alpha}}{n!\alpha^n}, \cdots \end{aligned} \tag{5.6}$$

Thus, by using (5.6) the approximate solution of (5.1) obtained by CADM is

$$\tilde{u}_m(x,t) = \sum_{n=0}^m u_n(x,t) = \sum_{n=0}^m (-1)^n \sin\left(\frac{x^\beta}{\beta}\right)\frac{t^{n\alpha}}{n!\alpha^n}. \tag{5.7}$$

From (5.7) we obtain

$$u(x,t) = \lim_{m\to\infty}\tilde{u}_m(x,t) = \sin\left(\frac{x^\beta}{\beta}\right)e^{-\frac{t^\alpha}{\alpha}}. \tag{5.8}$$

This analytical approximate solution (5.8) is the exact solution.

*Solution by CRDTM:*

Now solve this problem by using CRDTM. By taking the conformable fractional reduced differential transform (CRDT) of (4.1), it can be obtained that

$$\alpha(k+1)U_{k+1}^{\alpha}(x) = \frac{\partial^{2\beta}}{\partial x^{2\beta}}U_k^{\alpha}(x) \tag{5.9}$$

where the $t$-dimensional spectrum function $U_k^{\alpha}(x)$ is the conformable fractional reduced differential transform function. From the initial condition (5.2) we write

$$U_0^{\alpha}(x) = \sin\left(\frac{x^{\beta}}{\beta}\right) \tag{5.10}$$

Substituting (5.10) into (5.9), we obtain the following $U_k^{\alpha}(x)$ values successively

$$U_1^{\alpha}(x) = -\sin\left(\frac{x^{\beta}}{\beta}\right)\frac{1}{\alpha}, \quad U_2^{\alpha}(x) = \sin\left(\frac{x^{\beta}}{\beta}\right)\frac{1}{2!\alpha^2}, \cdots, U_n^{\alpha}(x) = \sin\left(\frac{x^{\beta}}{\beta}\right)\frac{(-1)^n}{n!\alpha^3}, \cdots$$

Then, the inverse transformation of the set of values $\{U_k^{\alpha}(x)\}_{k=0}^{n}$ gives the following approximation solution

$$\tilde{u}_n(x,t) = \sum_{k=0}^{n} U_k^{\alpha}(x)t^{k\alpha} = \sum_{k=0}^{n} \sin\left(\frac{x^{\beta}}{\beta}\right)\frac{(-1)^k}{k!\alpha^k}t^{k\alpha}. \tag{5.11}$$

From (5.11) we obtain

$$u(x,t) = \lim_{m\to\infty}\tilde{u}_m(x,t) = \sin\left(\frac{x^{\beta}}{\beta}\right)e^{-\frac{t^{\alpha}}{\alpha}}. \tag{5.12}$$

This analytical approximate solution (5.12) is the exact solution.
**Remark 4.1.**
If take $\alpha = \beta = 1$ in the problem, then Example 4.1 is reduced to standard diffusion equation

$$\frac{\partial}{\partial t}u(x,t) = \frac{\partial^2}{\partial x^2}u(x,t) \; t>0, x\in R$$

with initial condition

$$u(x,0) = \sin(x)$$

and our analytical approximate solution (5.8) and (5.12) imply

$$u(x,t) = \sin(x)e^{-t}$$

and this result is the exact solution of the standard problem in the literature.

The Aproximate solutions obtained with both CADM and CRDTM give us the exact solution.

**Example 4.2:**

Now let us consider the nonlinear time- and space fractional gas dynamics equation:

$$\frac{\partial^{\alpha}}{\partial t^{\alpha}}u(x,t) + \frac{1}{2}\frac{\partial^{\beta}}{\partial x^{\beta}}u^2(x,t) - u(x,t)(1-u(x,t)) = 0, \; 0<\alpha,\beta\leq 1 \tag{5.13}$$

subject to initial condition

$$u(x,0) = e^{-\frac{x^{\beta}}{\beta}}. \tag{5.14}$$

The exact solutions of (5.13) in conformable sense is

$$u(x,t) = e^{\frac{t^{\alpha}}{\alpha} - \frac{x^{\beta}}{\beta}}.$$

### Solution by CADM:

Solve the problem by using CADM. Let $L_\alpha = T_\alpha = \dfrac{\partial^\alpha}{\partial t^\alpha}$ be a linear operator, then the operator form of (5.13) is as follows

$$T_\alpha u(x,t) = -\frac{1}{2}\frac{\partial^\beta}{\partial x^\beta}u^2(x,t) + u(x,t)(1-u(x,t)), \quad 0<\alpha,\beta \leq 1. \tag{5.15}$$

By the help of Lemma 2.2, eq. (5.15) can be written as

$$t^{1-\alpha}\frac{\partial u(x,t)}{\partial t} = u(x,t) - u(x,t)\frac{\partial^\beta}{\partial x^\beta}u(x,t) - u^2(x,t), \quad 0<\alpha,\beta \leq 1. \tag{5.16}$$

If $L_\alpha^{-1} = \int_0^t \dfrac{1}{\xi^{1-\alpha}}(.)d\xi$, which is the inverse of $L_\alpha$, is applied to both sides of eq. (5.16), we get

$$u(x,t) = u(x,0) + L_\alpha^{-1}(u(x,t)) - L_\alpha^{-1}\left(u(x,t)\frac{\partial^\beta}{\partial x^\beta}u(x,t) + u^2(x,t)\right).$$

According to (3.8) and initial condition (5.14), we can write the following recursive relations:

$$\begin{aligned}
u_0 &= e^{-\frac{x^\beta}{\beta}} \\
u_1 &= L_\alpha^{-1}(u_0) - L_\alpha^{-1}(A_0) \\
&\vdots \\
u_{n+1} &= L_\alpha^{-1}(u_n) - L_\alpha^{-1}(A_n), n \geq 0.
\end{aligned} \tag{5.17}$$

where $A_n$'s are Adomian polynomials. By using the Maple code above, for the nonlinear term $N(u(x)) = u(x,t)\dfrac{\partial}{\partial x}u(x,t) + u^2(x,t)$, the Adomian polynomials can be obtain as:

$$\begin{aligned}
A_0 &= u_0^2 + u_0\frac{\partial^\beta}{\partial x^\beta}u_0 \\
A_1 &= 2u_0u_1 + u_0\frac{\partial^\beta}{\partial x^\beta}u_1 + u_1\frac{\partial^\beta}{\partial x^\beta}u_0 \\
A_2 &= u_1^2 + 2u_0u_2 + u_0\frac{\partial^\beta}{\partial x^\beta}u_2 + u_1\frac{\partial^\beta}{\partial x^\beta}u_1 + u_2\frac{\partial^\beta}{\partial x^\beta}u_0 \\
A_3 &= 2u_1u_2 + 2u_0u_3 + u_0\frac{\partial^\beta}{\partial x^\beta}u_3 + u_1\frac{\partial^\beta}{\partial x^\beta}u_2 + u_2\frac{\partial^\beta}{\partial x^\beta}u_1 + u_3\frac{\partial^\beta}{\partial x^\beta}u_0 \\
&\vdots
\end{aligned} \tag{5.18}$$

From (5.17) and (5.18), we conclude the terms of decomposition series as:

$$u_0 = e^{-\frac{x^\beta}{\beta}}, \; u_1 = e^{-\frac{x^\beta}{\beta}}\frac{t^\alpha}{\alpha}, \; u_2 = e^{-\frac{x^\beta}{\beta}}\frac{t^{2\alpha}}{2\alpha^2}, \; u_3 = e^{-\frac{x^\beta}{\beta}}\frac{t^{3\alpha}}{3!\alpha^3}, \cdots, u_n = e^{-\frac{x^\beta}{\beta}}\frac{t^{n\alpha}}{n!\alpha^n}, \cdots \tag{5.19}$$

Thus, from (5.19), the approximate solution of (5.13) obtained by CADM is

$$\tilde{u}_m(x,t) = \sum_{n=0}^m u_n(x,t) = \sum_{n=0}^m e^{-\frac{x^\beta}{\beta}}\frac{t^{n\alpha}}{n!\alpha^n}. \tag{5.20}$$

From (5.20) we obtain

$$u(x,t) = \lim_{m\to\infty}\tilde{u}_m(x,t) = e^{\frac{t^\alpha}{\alpha} - \frac{x^\beta}{\beta}}. \tag{5.21}$$

This analytical approximate solution (5.21) is the exact solution.

### Solution by CRDTM:

Now solve this problem by using CFDTM. By taking the conformable fractional reduced differential transform (CRDTM) of (5.13), it can be obtained that

$$\alpha(k+1)U_{k+1}^\alpha(x) = -\sum_{r=0}^k U_{k-r}^\alpha(x)\frac{\partial^\beta}{\partial x^\beta}U_r^\alpha(x) + U_k^\alpha(x) - \sum_{r=0}^k U_{k-r}^\alpha(x)U_r^\alpha(x) \tag{5.22}$$

where the $t$-dimensional spectrum function $U_k^\alpha(x)$ is the conformable fractional reduced differential transform function. From the initial condition (5.14) we write

$$U_0^\alpha(x) = e^{-\frac{x^\beta}{\beta}} \qquad (5.23)$$

Substituting (5.23) into (5.22), we obtain the following $U_k^\alpha(x)$ values successively

$$U_1^\alpha(x) = e^{-\frac{x^\beta}{\beta}}\frac{1}{\alpha},\ U_2^\alpha(x) = e^{-\frac{x^\beta}{\beta}}\frac{1}{2!\alpha^2},\cdots,U_n^\alpha(x) = e^{-\frac{x^\beta}{\beta}}\frac{1}{n!\alpha^3},\cdots$$

Then, the inverse transformation of the set of values $\{U_k^\alpha(x)\}_{k=0}^n$ gives the following approximation solution

$$\tilde{u}_n(x,t) = \sum_{k=0}^{n} U_k^\alpha(x) t^{k\alpha} = \sum_{k=0}^{n} e^{-\frac{x^\beta}{\beta}}\frac{1}{k!\alpha^k}t^{k\alpha}. \qquad (5.24)$$

From (5.24) we obtain

$$u(x,t) = \lim_{n\to\infty}\tilde{u}_n(x,t) = e^{\frac{t^\alpha}{\alpha}-\frac{x^\beta}{\beta}}. \qquad (5.25)$$

This analytical approximate solution (5.25) is the exact solution.

**Remark 4.2.**
If take $\alpha = \beta = 1$ in the problem, then Example 4.2 is reduced to standard gas dynamics equation

$$\frac{\partial}{\partial t}u(x,t) + \frac{1}{2}\frac{\partial}{\partial x}u^2(x,t) - u(x,t)(1-u(x,t)) = 0$$

with initial condition

$$u(x,0) = e^{-x}$$

our analytical approximate solution (5.21) and (5.25) imply

$$u(x,t) = e^{t-x}$$

and this result is the exact solution of the standard problem in the literature.

The approximate solutions obtained with both CADM and CRDTM give us the existed exact solution.

**Example 4.3:**

Now let us consider the nonlinear time- and space fractional partial differential equation

$$\frac{\partial^\alpha}{\partial t^\alpha}u(x,t) + (1+u(x,t))\frac{\partial^\alpha}{\partial x^\alpha}u(x,t) = 0,\ 0<\alpha\leq 1 \qquad (5.26)$$

subject to initial condition

$$u(x,0) = \frac{x^\alpha - \alpha}{2\alpha}. \qquad (5.27)$$

The exact solutions of (5.26) in conformable sense is

$$u(x,t) = \frac{x^\alpha - t^\alpha - \alpha}{t^\alpha - 2\alpha}.$$

***Solution by CADM:***

Solve the problem by using CADM. Let $L_\alpha = T_\alpha = \frac{\partial^\alpha}{\partial t^\alpha}$ be a linear operator, then the operator form of (5.26) is as follows

$$T_\alpha u(x,t) = -(1+u(x,t))\frac{\partial^\alpha}{\partial x^\alpha}u(x,t),\ 0<\alpha\leq 1. \qquad (5.28)$$

By the help of Lemma 2.2, eq. (5.28) can be written as

$$t^{1-\alpha}\frac{\partial u(x,t)}{\partial t} = -\frac{\partial^\alpha}{\partial x^\alpha}u(x,t) - u(x,t)\frac{\partial^\alpha}{\partial x^\alpha}u(x,t),\ 0<\alpha\leq 1. \qquad (5.29)$$

If $L_\alpha^{-1} = \int_0^t \frac{1}{\xi^{1-\alpha}}(.)d\xi$, which is the inverse of $L_\alpha$, is applied to both sides of eq. (5.29), we get

$$u(x,t) = u(x,0) - L_\alpha^{-1}\left(\frac{\partial^\alpha}{\partial x^\alpha}u(x,t)\right) - L_\alpha^{-1}\left(u(x,t)\frac{\partial^\alpha}{\partial x^\alpha}u(x,t)\right).$$

According to (3.8) and initial condition (5.27), we can write the following recursive relations:

$$u_0 = \frac{x^\alpha - \alpha}{2\alpha}$$
$$u_1 = L_\alpha^{-1}(u_0) - L_\alpha^{-1}(A_0) \qquad (5.30)$$
$$\vdots$$
$$u_{n+1} = L_\alpha^{-1}(u_n) - L_\alpha^{-1}(A_n), n \geq 0.$$

where $A_n$'s are Adomian polynomials. By using the Maple code above, for the nonlinear term $N(u(x)) = u(x,t)\frac{\partial}{\partial x}u(x,t) + u^2(x,t)$, the Adomian polynomials can be obtain as:

$$A_0 = u_0 \frac{\partial^\alpha}{\partial x^\alpha} u_0$$
$$A_1 = u_0 \frac{\partial^\alpha}{\partial x^\alpha} u_1 + u_1 \frac{\partial^\alpha}{\partial x^\alpha} u_0$$
$$A_2 = u_0 \frac{\partial^\alpha}{\partial x^\alpha} u_2 + u_1 \frac{\partial^\alpha}{\partial x^\alpha} u_1 + u_2 \frac{\partial^\alpha}{\partial x^\alpha} u_0 \qquad (5.31)$$
$$A_3 = u_0 \frac{\partial^\alpha}{\partial x^\alpha} u_3 + u_1 \frac{\partial^\alpha}{\partial x^\alpha} u_2 + u_2 \frac{\partial^\alpha}{\partial x^\alpha} u_1 + u_3 \frac{\partial^\alpha}{\partial x^\alpha} u_0$$
$$\vdots$$

From (5.30) and (5.31), we conclude the terms of decomposition series as:

$$u_0 = \frac{x^\alpha - \alpha}{2\alpha},\ u_1 = -\frac{x^\alpha + \alpha}{(2\alpha)^2}t^\alpha,\ u_2 = \frac{x^\alpha + \alpha}{(2\alpha)^3}t^{2\alpha},$$
$$u_3 = -\frac{x^\alpha + \alpha}{(2\alpha)^4}t^{3\alpha}, \cdots, u_n = (-1)^n \frac{x^\alpha + \alpha}{(2\alpha)^{n+1}}t^{n\alpha}, \cdots \qquad (5.32)$$

Thus, from (5.32) the approximate solution of (5.26) obtained by CADM is

$$\tilde{u}_m(x,t) = \sum_{n=0}^{m} u_n(x,t) = \frac{x^\alpha - \alpha}{2\alpha} + \sum_{n=1}^{m}(-1)^m \frac{x^\alpha + \alpha}{(2\alpha)^{m+1}}t^{m\alpha}. \qquad (5.33)$$

*Solution by CRDTM:*

Now solve this problem by using CFDTM. By taking the conformable fractional reduced differential transform (CRDTM) of (5.26), it can be obtained that

$$\alpha(k+1)U_{k+1}(x) = -U_k(x) - \sum_{r=0}^{k} U_{k-r}(x)\frac{\partial}{\partial x}U_r(x) \qquad (5.34)$$

where the $t$-dimensional spectrum function $U_k^\alpha(x)$ is the conformable fractional reduced differential transform function. From the initial condition (5.27) we write

$$U_0^\alpha(x) = \frac{x^\alpha - \alpha}{2\alpha} \qquad (5.35)$$

Substituting (5.35) into (5.34), we obtain the following $U_k^\alpha(x)$ values successively

$$U_1^\alpha(x) = -\frac{x^\alpha + \alpha}{(2\alpha)^2},\ U_2^\alpha(x) = \frac{x^\alpha + \alpha}{(2\alpha)^3}, \cdots, U_n^\alpha(x) = (-1)^n \frac{x^\alpha + \alpha}{(2\alpha)^{n+1}}, \cdots$$

Then, the inverse transformation of the set of values $\{U_k^\alpha(x)\}_{k=0}^{n}$ gives the following approximation solution

$$\tilde{u}_m(x,t) = \sum_{k=0}^{m} U_k^\alpha(x)t^{k\alpha} = \frac{x^\alpha - \alpha}{2\alpha} + \sum_{n=1}^{m}(-1)^m \frac{x^\alpha + \alpha}{(2\alpha)^{m+1}}t^{m\alpha}. \qquad (5.36)$$

Compare the seventh iteration CADM and CFDTM solutions with the exact solution on the graph for some $\alpha$ values.

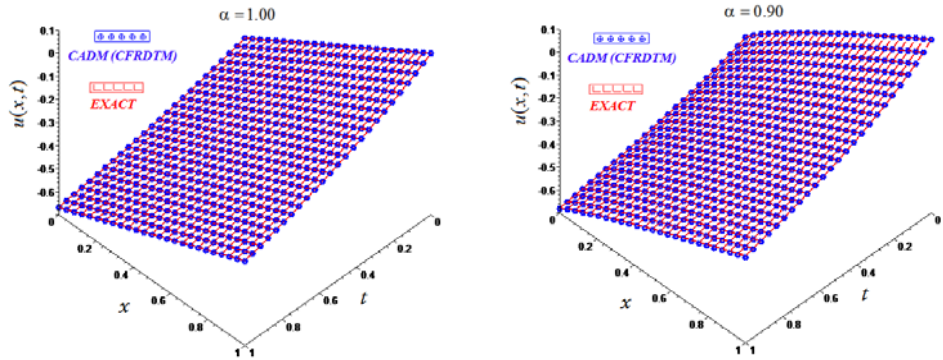

**Figure 6.1.** Comparison of the seventh iteration approximate solutions of CADM (CRDTM) with the exact solutions for eq. (5.26).

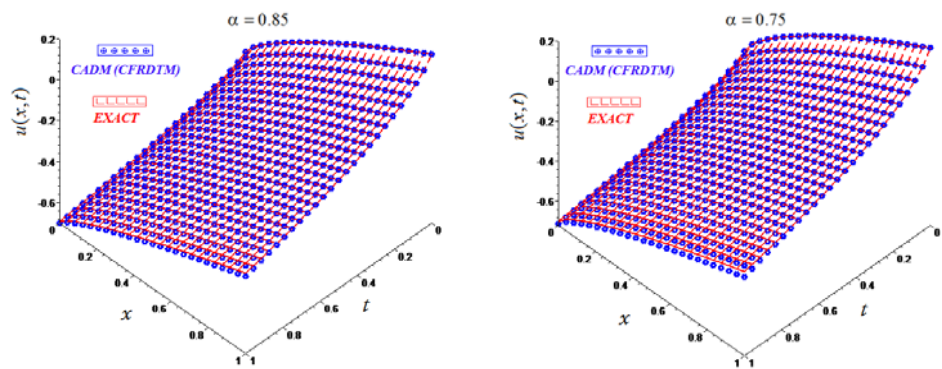

**Figure 6.2.** Comparison of the seventh iteration approximate solutions of CADM (CRDTM) with the exact solutions for eq. (5.26)

## 6. Conclusion

The fundamental goal of this article is to construct the approximate solutions of FPDEs. The goal has been achieved by using CADM for the first time and it is compared with CRDTM. CADM and CRDTM are applied to different linear and non-linear conformable time and space FPDEs. And also the approximate analytical solutions obtained by CADM and CRDTM are compared to each other and with the exact solutions. CADM and CRDTM offer solutions with easily computable components as convergent series. Approximate solutions obtained by CADM are exactly same as the solutions obtained by CRDTM for time and space FPDEs. The CADM gives quantitatively reliable results as CRDTM, and also it requires less computational work than existing other methods. As a result, in recent years, FDEs emerging as models in fields such as mathematics, physics, chemistry, biology and engineering makes it necessary to investigate the methods of solutions and we hope that this study is an improvement in this direction.